
\documentstyle{amsppt}
\magnification 1200
\NoBlackBoxes
\NoRunningHeads

\define\a{\alpha}
\define\la{\lambda}
\define\La{\Lambda}
\define\ver{M_\lambda}
\define\op#1{\operatorname{#1}}
\define\qq#1{{{\bmatrix #1 \endbmatrix}_q}}

\define\tw{\Phi_\lambda^u}
\define\twin{\tilde\Phi_\lambda^u}
\define\End{\op{End\big(U[0]\big)}}

\define\hw{v_\lambda}
\define\lw{v_{-\lambda}^*}

\redefine\o{\otimes}
\define\({\left(}
\define\){\right)}
\define\[{\left[}
\define\]{\right]}
\define\<{\langle}
\define\>{\rangle}
\define\qgr{\Cal U_q\frak g}
\define\g{\frak g}
\define\h{\frak h}
\define\Z{{\Bbb Z}}
\define\N{{\Bbb N}}

\define\C{{\Bbb C}}
\redefine\P{{\bold P}}
\define\Q{{\bold Q}}

\topmatter
\title Algebraic integrability of Macdonald operators
       and representations of quantum groups 
  \endtitle

\author
Pavel Etingof
\footnote{Department of Mathematics, Harvard University,
          Cambridge, MA 02138, USA; \newline
          e-mail: etingof\@math.harvard.edu \newline}
and Konstantin Styrkas
\footnote{Department of Mathematics, Yale University,
          New Haven, CT 06520, USA; \newline
          e-mail: styrkas\@math.yale.edu \newline}
\endauthor

\abstract
In this paper we construct examples of commutative rings of difference
operators with matrix coefficients from representation theory
of quantum groups, generalizing the results of our previous paper
\cite{ES} to the $q$-deformed case. A generalized Baker-Akhiezer 
function $\Psi$ is realized as a matrix character of a Verma module
and is a common eigenfunction for a commutative ring
of difference operators.

In particular, we obtain the following result in Macdonald theory:
at integer values of the Macdonald parameter $k$, there exist
difference operators commuting with Macdonald operators 
which are not polynomials of Macdonald operators. 
This result generalizes an analogous result of Chalyh and
Veselov for the case $q=1$, to arbitrary $q$.  
As a by-product, we prove a generalized Weyl character formula for
Macdonald polynomials (= Conjecture 8.2 from \cite{FV}), 
the duality for the $\Psi$-function, and
 the existence of shift operators.
\endabstract

\endtopmatter

\document
\vskip .1in
\centerline{March 25, 1996}
\vskip .1in
\centerline{\bf 1. Introduction.}
\vskip .1in

Let $N$ be a positive integer. Let $\frak D_q^N$ be the algebra over the
field $\C(q)$ generated by the field of rational functions
$\C(q,X_1,...,X_N)$ and commuting operators $T_1^{\pm 1},...,T_N^{\pm 1}$,
with commutation relations 
$$T_i\circ f(q,X_1,...X_i,...,X_N)=f(q,X_1,...,qX_i,...,X_N)\circ T_i.$$ 
This algebra is called the algebra of $q$-difference operators 
in $n$ variables with rational coefficients.
Elements of this algebra are called difference operators. 

Let $V$ be a finite-dimensional vector space over $\C$. Introduce the algebra
$\frak D_q^N(V)$ of difference operators with matrix coefficients 
$$\frak D_q^N(V) = \frak D_q^N\o\op{End(V)}$$

Let $\g$ be a simple finite-dimensional Lie algebra over $\C$ of rank
$r$, and let $\qgr$ be the corresponding quantum group. In \cite{EK}, to any
finite dimensional representation $U$ of $\qgr$ was assigned
a family of commuting difference operators $D_c$ parametrized
by Weyl group invariant trigonometric polynomials $c$ on the Cartan
subalgebra $\h$ of $\g$.
These operators are constructed as follows.

Let $M_\la$ be the Verma module over $\qgr$ with highest weight
$\la$ and highest weight vector $v_\la$.
Let $U[0]$ be the zero weight subspace of $U$. For any $u\in U[0]$,
define the intertwining operator $\Phi_\la^u: M_\la\to M_\la\o U$ by the
condition $\Phi_\la^uv_\la=v_\la\o u+\sum w_i\o u_i$, where 
$w_i$ are homogeneous vectors of weights $\mu_i<\la$.
This operator is defined for generic $\la$.
For any weight $\nu$, let $\op{Proj}\bigl|_{\ver[\nu]}: \ver\to\ver$ be the
homogeneous projector to the subspace $\ver[\nu]$ of weight $\nu.$
Let $\tilde\psi_\la(X_1,...,X_r)$ be the function with values
in $\End$ such that for any $u\in U[0]$
$$\tilde\psi_\la(X_1,...,X_r)u=\sum_\mu X_1^{\mu_1}...X_r^{\mu_r}
Tr\bigl|_{M_\la}(\op{Proj}\bigl|_{\ver[\mu]} \circ \Phi_\la^u \circ
\op{Proj}\bigl|_{\ver[\mu]})$$

\proclaim{Proposition 1.1} \cite{EK} For any Weyl group invariant   
function $c(\la)$ on $\h^*$ of the form 
$$c(\la)=\sum_{\mu\in \P}c_\mu q^{2\<\la,\mu\>},\quad c_\mu\in\C(q)\tag 1-1 $$
there exists a unique difference operator
$D_c\in \frak D_q^r(U[0])$ such that
$$D_c \tilde\psi_\la = c(\la+\rho)\tilde\psi_\la$$
\endproclaim

For any root $\alpha$ of $\g$, let $k_\alpha=\max\{n|U[n\a]\ne0\},$
where $U[\mu]$ is the subspace of weight $\mu$ in $U$. 
Let $R(U)$ be the ring of functions on $\h^*$ of the form (1-1)
such that for any positive root $\alpha$ of $\g$
$$c\(\la-\frac{n\alpha}2 \)=c\(\la+\frac{n\alpha}2 \),\qquad 
n=1,...,k_\alpha $$
whenever $\<\la,\alpha\>=0.$
The main result of this paper is the following theorem, proved in 
Chapter 4 of this paper. 

\proclaim{Theorem 1.2}
There exists an injective homomorphism $\xi: R(U)\to \frak D_q^r(U[0])$
such that for any Weyl group invariant element $c\in R$ one has $\xi(c)=D_c$.
For any $c\in R(U)$, the operator $\xi(c)$ is defined by the equation
$$\xi(c)\tilde\psi_\la=c(\la+\rho)\tilde\psi_\la.$$
\endproclaim

We will denote $\xi(c)$ by $D_c$ for any $\xi\in R$.

In the case when $\g=\frak{sl}_N$ (type $A_{N-1}$), we can
choose representation $U$ to be $S^{kN}V$, where $V$ is the fundamental
representation, in which case
the space $U[0]$ is 1-dimensional.  Then the operators
$D_c$ for symmetric functions $c$  are conjugate
to Macdonald operators, corresponding to $t=q^{k+1}.$
Namely, if $c_l$ are elementary symmetric
functions, then $\{D_{c_l}\}$ are simultaneously conjugate to
$$\Cal M_l=\sum_{I\subset {1,...,N}, |I|=l}\prod_{i\in I,j\notin I}
\frac{q^{k+1}X_i-q^{-k-1}X_j}{X_i-X_j}\prod_{i\in I}T_i^2$$
(in suitable coordinates). In this case, the numbers $k_\alpha$ are
all equal to $k$; so we will denote the algebra $R(U)$ by $R_k$.
From Theorem 1.2 we get (See Chapter 5):
\proclaim{Theorem 1.3} For any positive integer $k$,
there exists an injective homomorphism $\xi:R_k\to \frak D_q^N$ such that 
$\xi(c_l)=\Cal M_l$, $l=1,\dots,N$. The function $\tilde\psi_\la$ is a common
eigenfunction of the operators $\xi(c), c\in R_k$ with eigenvalue
$c(\la+\rho)$. 
\endproclaim

Note that Theorem 1.3 is a special property of Macdonald's operators
at integer values of $k$. If $k$ is not an integer, one can show that
the centralizer of $\Cal M_1,\dots,\Cal M_N$ in $\frak D_q^N$
reduces to the polynomial algebra of $\Cal M_1,\dots,\Cal M_N$.
We call this special property at integer values of $k$ ``algebraic 
integrability of Macdonald operators'', by analogy with the
case differential operators which was treated in
\cite{CV1,CV2,VSC,ES}. In this sense, the results of this paper is 
precisely a $q$-deformation of the results of \cite{ES}.

As a by-product, we obtain several results in
Macdonald's theory. Namely, we prove the partial Weyl group
symmetry of the $\tilde\psi$-function, a generalized Weyl character
formula for Macdonald's polynomials
(which coincider with Conjecture 8.2 in
\cite{FV}), an explicit formula for the
$\tilde\psi$-function in terms of shift operators, and symmetry of the
$\tilde\psi$-function with respect to the interchange $\la\leftrightarrow x$.

The paper is organized as follows. In Section 2 we recall
basic facts about representations of quantum groups and
intertwining operators. In Section 3 we introduce the 
$\Psi$-function as matrix trace of an intertwining operator,
and prove its properties. In Section 4 we explain how to
construct a commutative ring of difference operators from
the $\Psi$-function. In Section 5 we review some facts from Macdonald
theory for root system $A_n$ and explain how to obtain them
from our construction. In Appendix we show how our technique
works in the simplest example.
\vskip .1in

\centerline{\bf Acknowledgements}

We would like to thank Igor Frenkel and Alexander Kirillov Jr. for
useful discussions. The work of the first author was partially supported by
an NSF postdoctoral fellowship.

\vskip .1in
\vskip .1in
\centerline{\bf 2. Quantum groups and their representations}
\vskip .1in

\subhead Notation \endsubhead
Let $\frak g$ be a simple (finite-dimensional) complex Lie algebra
of rank $r$ with fixed diagonalizable Cartan matrix
$A=(a_{ij}), i,j=1,\dots,r,$ and let $d_1,\dots,d_r$ be positive relatively
prime integers such that the matrix $B = (b_{ij}) =(d_ia_{ij})$
is symmetric. We denote its Cartan subalgebra by $\h.$
Let $\a_i\in\h^*, i=1,\dots,r$ denote simple roots,
$R$ - the corresponding root system,
$R^+$ and $R^-$ - the sets of positive and negative roots, respectively.

The invariant Killing form $\<\cdot,\cdot\>$ on $\h^*$ is defined by
$\<\a_i,\a_j\> = d_i a_{ij}.$
Let $\La_1,\dots,\La_r\in \h^*$ be fundamental weights, i.e.
$\<\La_i,\a_j\> = \delta_{ij},\quad i,j = 1,\dots,r.$
Put $\rho = \sum_{i=1}^r \La_i.$ Denote
$$\Q = \sum \Z\a_i,\Q_+ = \sum \Z_+\a_i,
\P = \sum \Z\La_i, \P_+ = \sum \Z_+\La_i.$$
For $\mu,\nu\in \P$ we write $\mu\ge\nu$ if $\mu-\nu\in \Q_+.$

Let $W$ be the Weyl group of $\frak g$. The
Weyl group generators $s_i$ act on $\h^*$ by simple root reflections
$$s_i\cdot\mu = \mu - 2\frac{\<\a_i,\mu\>}{\<\alpha_i,\alpha_i\>}\a_i$$
We also introduce a shifted action of Weyl group by
$$w^\rho \cdot \mu = w(\mu+\rho) - \rho$$
For $w\in W$ let $l(w)$ denote the length of $w,$
i.e. the number of generators in a reduced
decomposition $w = s_{i_1} \cdot \dots \cdot s_{i_l}.$

\subhead Quantum groups \endsubhead
The quantum group $\qgr,$ associated to simple Lie algebra $\g,$ is a
Hopf algebra over $\C(q)$ with generators
$E_i,F_i,K_i,$ $i=1,\dots,r$ and relations:
$$K_iK_j=K_jK_i,\quad K_iE_j=q_i^{a_{ij}}E_jK_i,\quad K_iF_j=q_i^{-a_{ij}}F_jK_i$$
$$E_iF_j-F_jE_i=\delta_{ij}\frac{K_i-K_i^{-1}}{q_i-q_i^{-1}},$$
$$\sum_{k=0}^{1-a_{ij}} (-1)^k \bmatrix 1-a_{ij} \\ k \endbmatrix_{q_i}
E_i^{1-a_{ij}-k}E_jE_i^k=0,\quad i\ne j,$$
$$\sum_{k=0}^{1-a_{ij}} (-1)^k \bmatrix 1-a_{ij} \\ k \endbmatrix_{q_i}
F_i^{1-a_{ij}-k}F_jF_i^k=0,\quad i\ne j.$$
where $q_i = q^{d_i}$ and we used notation
$$ \bmatrix n \\ k \endbmatrix _q=\frac{[n]_q!}{[k]_q! [n-k]_q!},
\quad  [n]_q! = [1]_q \cdot [2]_q \cdot \dots \cdot [n]_q,
\quad  [n]_q= \frac {q^n - q^{-n}}{q-q^{-1}}$$

Comultiplication $\Delta,$ antipode $S$ and counit $\epsilon$
in $\qgr$ are given by
$$\Delta(E_i) = E_i\otimes 1 + K_i\otimes E_i, \quad
  \Delta(F_i) = F_i\otimes K_i^{-1} + 1\otimes F_i, \quad
  \Delta(K_i) = K_i \otimes K_i$$
$$S(E_i)=-K_i^{-1}E_i,\quad S(F_i)=-F_iK_i,\quad S(K_i)=K_i^{-1}$$
$$\epsilon(E_i) = \epsilon(F_i) = 0,\quad \epsilon(K_i) = 1$$
We define a $\C$-algebra involution $\omega$ of $\qgr$ by
$$\omega(E_i)=-F_i,\quad\omega(F_i)=-E_i,\quad\omega(K_i)=K_i,
\quad\omega(q) = q^{-1}$$
We have a decomposition of vector spaces
$\qgr = \Cal U^-\otimes \Cal U^0 \otimes \Cal U^+,$
where $\Cal U^-$ (resp. $\Cal U^+$) is the subalgebra generated by $F_i$ (resp.
$E_i$), and $\Cal U^0$ is generated by $K_i,K_i^{-1}, i=1,\dots,r.$

\subhead Verma modules \endsubhead
For any $\la\in\h^*$ we can introduce Verma module $\ver$ over $\qgr,$
i.e. $\Cal U^-$-free module with a single generator $\hw$ and relations
$$E_i\hw=0,\quad K_i\hw=q^{\<\a_i,\la\>}\hw.$$

{\bf Remark.} Here and below we work over the field 
$F=\Bbb C(q)(\{q^a,a\in\Bbb C\})$. In this setting, $q^{\<\mu,\la\>}$ is
a function $\h^*\to F$. 
 
We have the decomposition
$$\ver=\bigoplus_{\mu\in \Q_+} \ver[\la-\mu]$$
of $\ver$ into direct sum of weight subspaces $\ver[\la-\mu],$
where we say that a vector $v$ has weight $\mu\in\h^*$ if
$$K_iv=q^{\<\a_i,\mu\>}v$$

The restricted dual module $\ver^*$ is a $\Cal U^+$-module with a
lowest weight vector $\lw$ such that $\<\lw,\hw\>=1.$ By definition we have
$$\<gv^*,v\> = \<v^*,S(g)v\>,\quad v\in\ver, v^*\in\ver^*.$$

Introduce a symmetric form
$F_\la$ on $\ver$ defined by
$$F(g_1\hw,g_2\hw)=\<\omega(g_1)\lw,g_2\hw\>.$$
The weight subspaces are pairwise orthogonal with respect
to this form. The restriction of $F$ to weight subspaces
$M_{\la-\mu}$ is proportional to the quantum Shapovalov 
form $\Cal F,$ introduced in \cite{CK}:
$$F_\mu(\cdot,\cdot) = C_\mu q^{-\<\la,\mu\>} \Cal F_\mu(\cdot,\cdot)$$
for some constants $C_\mu.$

Fix a basis $g_i^\mu\in\Cal U^-[\mu].$
Let $F_\mu = (F_\mu)_{ij},\quad i,j=1,2,\dots,\dim\ver[\la-\mu],$ denote
the matrix of the restriction of form $F$ to $\ver[\la-\mu]$ with
respect to the basis $g^\mu_i \hw$ of $\ver[\la-\mu].$
A variation of the quantum determinant formula [CK]
asserts that
$$\det F_\mu(\la) = C \prod_{\a\in R^+}\prod_{n\in\N}
\( 1 - q^{-2\<\a,\la+\rho\> + n\<\a,\a\>} \)^{\op{Par}(\mu-n\a)},$$
where $\op{Par}$ is the generalized Kostant partition function,
and $C\in\C(q)$ is a constant, depending on the choice of basis $g_i^\mu.$

This determinant is a linear combination of terms
$q^{-2\<\mu,\la\>},$ where $\mu$'s belong to a finite
subset $L \subset \Q,$ with some coefficients from $\Bbb C(q).$
Motivated by this fact,
we introduce

\proclaim{Definition}
Expressions of the form $\sum_{\mu\in L} a_\mu q^{-2\<\mu,\la\>}$,
$a_\mu\in\Bbb C(q)$,
will be called $q$-polynomials with support $L$ and coefficients
$a_\mu.$
\endproclaim

Verma modules are reducible when form $F$ is degenerate,
i.e. $\det F_\mu(\la) = 0$ for some $\mu.$
This happens when $\la$ satisfies one of the Kac-Kazhdan equations:
$$\<\a,\la+\rho\>=\frac n2\<\a,\a\>,\quad n=1,2,\dots \tag 2-1$$
For $\la$ generic from Kac-Kazhdan hyperplanes,
$\ver$ contains a unique submodule $\ver^1,$ isomorphic to
$M_{\la-n\a}.$

\subhead Intertwining operators \endsubhead
Let $U$ be an irreducible finite-dimensional $\qgr$-module with
non-trivial zero weight subspace $U[0].$ For $u\in U$ let
$\tw:\ver\to\ver\otimes U$ be an intertwining operator such that
$\hw\to\hw\otimes u + \text{higher order terms,}$
where ``higher order terms'' mean terms of the form 
$v_{\la-\mu}\otimes u_\mu, \mu>0.$

If $\ver$ is irreducible, then $\tw$ exists and is unique for
any $u\in U[0].$ Indeed, we have a unique $\Cal U^+$-intertwiner
$\Omega: \ver^* \to U,$ such that $\Omega \lw = u.$ Since
$Hom(\ver^*,U) \cong \ver^{**} \otimes U \cong \ver\otimes U,$
it corresponds to a singular (i.e. $\Cal U^+$-invariant) vector
$\phi\in \ver\otimes U.$ We now construct $\tw$ by putting
$\tw\hw = \phi$ and extending $\tw$ to the whole $\ver$
by the intertwining property.

For our purposes we need an explicit form for that singular vector.

\proclaim{Proposition 2.1}
For any (homogeneous) basis $\{g^\mu_i\}$ of $\Cal U^-$
$$\phi =\sum_\mu \(\sum_{i,j} \(F_\mu^{-1}\)_{ij}(\la)
g_i^\mu\hw\otimes \omega(g_j^\mu)u\) \tag 2-2$$
is a singular vector in $\ver\otimes U.$
\endproclaim
Note that since $U$ has a highest weight, the summation is
over the finite set of $\mu$'s such that $U[\mu]\ne 0.$
\demo{Proof}
We check that the corresponding element $\Phi \in Hom(\ver^*,U)$
defined as the composition $\ver^* \to \ver^*\otimes\ver\otimes U\to U$
is a $\Cal U^+$-intertwiner. We have:
$$\gather
\Phi \(\omega(g^\nu_n \lw)\) = \sum_\mu \(\sum_{ij} \(F_\mu^{-1}\)_{ij}
 \<\omega(g^\nu_n)\lw,g^\mu_i\hw\> \omega(g^\mu_j)u \) = \\
 \sum_{i,j} \(F_\nu^{-1}\)_{ij}
 \<\omega(g^\nu_n)\lw,g^\nu_i\hw\> \omega(g^\nu_j)u =
 \sum_j \(\sum_i \(F_\nu^{-1}\)_{ij} \(F_\nu\)_{ni}\) \omega(g^\nu_j)u =\\
 \sum_j \delta_{jn} \omega(g^\nu_j)u = \omega(g^\nu_n) u
\endgather$$
\qed
\enddemo

Recall that in classical case matrix elements of the inverse matrix
$F_\mu^{-1}$ were rational functions of $\la$ with at most simple
poles in the Kac-Kazhdan hyperplanes given by (2-1). (see \cite{ES}).
A similar argument, also involving Jantzen filtration, proves that the
same is true in the quantum case. Therefore, expression (2-2) for
the singular vector $\phi$ is a ratio of two $q$-polynomials, with
at most simple singularities in a finite collection of
Kac-Kazhdan hyperplanes.

If we multiply the $q$-rational expression (2-2) by the least common
denominator $\tilde\chi(\la),$ we will get a well-defined for all
$\la$'s formula for a singular vector $\tilde\phi\in\ver\otimes U.$
We are now going to show that in fact the least common denominator may
only contain factors
$$\chi_n^\a(\la) = 1 - q^{-2\<\a,\la+\rho\>+n\<\a,\a\>},$$
corresponding to $n,\a$ such that $U[n\a]\ne0.$ Indeed,
suppose that $\tilde\chi(\la)$ contained a factor 
$\chi_n^\a(\la),$ but $U[n\a]=0.$

Consider $\la$ generic from the hyperplane
$\<\a,\la+\rho\> = \frac n2 \<\a,\a\>.$ Then $\ver$ contains a
unique maximal submodule $\ver^1 \cong M_{\la-n\a},$ generated
by singular vector $v_{\la-n\a}.$ Since the first term 
$\tilde\chi(\la) \hw\otimes u$ in the expression for $\tilde\phi$
turns into zero in our hyperplane, the singular vector must
have form
$$\tilde\phi = v_{\la-n\a} \otimes \tilde u + \text{higher order terms}$$
The intertwining property implies that $\tilde u \in U[n\a],$
and by assumption $\tilde u = 0.$ Therefore, $\tilde\phi$ is
zero for $\la$ generic from the hyperplane, and by Bezout theorem
is divisible by $\chi_n^\a.$ This shows that $\tilde\chi(\la)$
was not the least common denominator - contradiction.

Denote $k_\a = \max\{n|U[n\a]\ne0\},$
$$L_\theta = \{\nu|\nu = \sum_{\a\in R^+} m_\a \a, \quad 0\le m_\a \le k_\a
 \text{ for all } \a\in R^+ \} $$

We conclude this section with the following
\proclaim{Proposition 2.2}
If $U$ is an irreducible finite-dimensional $\qgr$-module with highest
weight $\theta,$ the singular vector 
$\phi\in\ver\otimes U,$ given by (2-2), can be represented as
$$\phi = \frac {\sum_l S_l(\la) \tilde g_l \hw \otimes u_l}
{\prod_{\a\in R^+}\prod_{m=1}^{k_\a}
  \(1 - q^{-2\<\a,\la+\rho\>+m\<\a,\a\>} \)}$$
where $g_l\in\Cal U^-, u_l\in U,$
and $q$-polynomials $S_l(\la)$ have supports, contained
in $L_\theta.$
\endproclaim
\demo{Proof}
We already proved that the least common denominator for the expression
for $\phi$ may only contain factors $\chi_n^\a(\la), n=1,\dots,k_\a.$

The statement about the support of the polynomials $S_l(\la)$ follows from
the fact that the support of the numerator must lie within the convex
hull of the support of the denominator, which in this case is exactly
$L_\theta.$
\qed\enddemo

\vskip .1in

\vskip .1in
\centerline{\bf 3. Matrix Trace, the $\Psi$-function and its properties.}
\vskip .1in

We now fix an irreducible $\qgr$-module $U$ with highest weight $\theta$ and
non-trivial zero weight subspace. We use the notation
$$k_\a = \max\{n | U[n\a] \ne 0\}, \qquad 
 \Theta = \sum_{\a\in R^+}k_\a\cdot\a \in\Q_+$$
$$ \chi_n^\a(\la) = 1 - q^{-2\<\a,\la+\rho\>+n\<\a,\a\>}, \qquad
 \chi(\la) = \prod_{\a\in R^+}\prod_{n=1}^{k_\a} \chi_n^\a(\la)$$
$$L_\theta = \{ \mu\in\Q_+ | \mu = \sum_{\a\in R^+} m_\a \a,
   \quad 0 \le m_\a \le k_\a \text{ for all } \a\in R^+ \}$$

As in [ES], define a new intertwining operator
$$\twin = \chi(\la) \tw : \ver \to \ver\otimes U $$
From Proposition 2.2 it follows that $\twin$ is
well-defined even when $\la$ belongs to Kac-Kazhdan hyperplanes,
where $\tw$ did not always exist.

Introduce a $\End$-valued function $\Psi(\la,x), \la,x\in\h^*,$ by
$$\Psi(\la,x) u = Tr\bigr|_{M_\la}(\twin e^x) $$
\proclaim {Proposition 3.1}
The function $\Psi(\la,x)$ defined above has the form
$$\Psi(\la,x) = e^{\<\la,x\>}
\sum_{\mu\in L_\theta}q^{-2\<\mu,\la+\rho\>} P_\mu(x),$$
where $P_\mu(x) \in \End$ and $P_\Theta(x)$ is invertible
for generic $x.$

If $\<\a,\la+\rho\> = \frac n2\<\a,\a\>$ for some 
$\a\in R^+, n=1,2,\dots,k_\a$ then
$$\Psi(\la,x) = \Psi(\la-n\a,x) B_{n\a}(\la) \tag 3-1$$
for some (possibly infinite) sum
$$B_{n\a}(\la)=\sum_{\mu\in \Q_+}q^{-2\<\mu,\la\>} B_{n\a}^\mu,
\qquad  B_{n\a}^\mu\in \End$$
(In fact, the matrix elements of
$B_{n\a}(\la)$'s are ratios of $q$-polynomials.)
\endproclaim

{\bf Remark.}
If we take $U$ to be a trivial module, then the $\Psi$-function
becomes the usual character of the Verma module. Therefore,
we can regard the $\Psi$-function as a generalized (matrix-valued)
character of the Verma module $\ver.$

\demo{Proof of Proposition 3.1}
$$Tr\bigr|_{\ver}\(\twin e^x\) = e^{\<\la,x\>}
\sum_{\mu\in\Q_+} e^{-\<\mu,x\>} \Cal B_\mu(\la)u,$$
where ``partial traces'' $\Cal B_\mu(\la)\in\End,$ corresponding
to weight subspaces $\ver[\la-\mu],$ are defined by
$$\Cal B_\mu(\la)u = Tr \(
 \op{Proj}_{\ver[\la-\mu]} \circ \twin
 \circ \op{Proj}_{\ver[\la-\mu]} \) $$
By Proposition 2.2, $\Cal B_\mu(\la)$ are $\End$-valued
$q$-polynomials with support $L_\theta.$ If we let $x\to\infty$ in
such a way that  $\<\a,x\>\to +\infty,\alpha\in R_+$ 
(this just means that we are keeping only the highest weight terms 
of the series), we will get asymptotically
$$\Psi(\la,x) \sim e^{\<\la,x\>} \chi(\la) \cdot \bold 1,$$
Therefore $P_\Theta(x) \sim \bold 1,$ and $P_\Theta(x)$ is
invertible for generic $x.$

We now prove the second property.

If $\la$ is generic from hyperplane $\<\a,\la+\rho\> = \frac n2\<\a,\a\>,$
then $\ver$ is reducible and contains a unique submodule $\ver^1$
generated by singular vector $v_{\la-n\a}.$ It is clear that for such
$\la$ there are no order $\mu$ terms in $\tw\hw$ unless $\mu\ge n\a.$
In other words, $\twin$ maps $\ver$ into $\ver^1\otimes U,$ and
$$\twin\hw = v_{\la-n\a} \otimes u' + \text{higher order terms}$$
$$\twin v_{\la-n\a} = v_{\la-n\a} \otimes u'' + \text{higher order terms}$$
Clearly, both $u'$ and $u''$ depend linearly on $u.$
More precisely,
$$u'' = \Cal B_{n\a}(\la) u, \tag 3-3$$
Therefore we have
$$\chi(\la-n\a) \twin\bigr|_{\ver^1} \cong \tilde\Phi_{\la-n\a}^{u''}$$
Taking the traces of the operators from the last equation and
using (3-3), we get
$$\chi(\la-n\a)\Psi(\la,x) = \Psi(\la-n\a,x) \Cal B_{n\a}(\la)$$
Since $\chi(\la-n\a)$ is invertible in the ``Laurent series''
completion of $\C[\P],$ we can introduce
$$B_{n\a}(\la) = \frac{\Cal B_{n\a}(\la)}{\chi(\la-n\a)} =
\sum_{\mu\in \Q_+} q^{2\<\mu,\la\>} B_{n\a}^\mu$$
and (3-1) follows. Proposition 3.1 is proved.
\qed\enddemo

We prove that property (3-1) of the $\Psi$-function
determines it uniquely up to multiplication by a factor,
depending only on $x.$

\proclaim {Proposition 3.2}
Suppose we have an $\End$-valued function
$$\Psi'(\la,x) = e^{\<\la,x\>} \sum_{\mu\in L} q^{-2\<\mu,\la\>} Q_\mu(x)$$
where $L \subset \Q_+$ and $Q_\mu(x)\in\End,$ satisfying condition
(3-1). Then the $L$ contains at least one weight $\mu\ge\Theta.$
\endproclaim

\demo{Proof (cf. [ES])}
Let us rewrite the condition (3-1). We have:
$$e^{\frac n2\<\a,x\>} \sum_{\mu\in \Q}
q^{-2\<\mu,\la+\frac{n\a}2\>} P_\mu(x)=
e^{-\frac n2\<\a,x\>} \sum_{\mu\in \Q}
q^{-2\<\mu,\la-\frac{n\a}2\>} P_\mu(x)
\sum_{\nu\in \Q}q^{-2\<\nu,\la\>} B_{n\a}^\nu.$$

For every $\mu\in \Q$ consider the set
$$Z_\a(\mu) = \mu +\Z\a = \{\nu\in \Q | \nu = \mu + m\a
\quad\text{ for some }m\in\Z \}$$

Comparing coefficients for
$q^{\<\mu,\la\>},$ we get for $n=1,2,\dots,k_\a$:
$$\sum_{\nu\in Z_\a(\mu)}
\( e^{n\<\a,x\>}q^{-n\<\nu,\a\>} P_\nu(x) - \sum_{\beta\in \Q}
q^{-n\<\beta,\a\>} P_{\nu-\beta}(x) B_{n\a}^\beta \) = 0 $$

This is a system of linear equations on unknown functions $P_\nu(x).$
Note that the summation over $\beta$ is finite, since $P_{\nu-\beta}=0$
if $\nu-\beta\notin \Q_+.$
The matrix of this system has a block-upper-triangular form,
blocks corresponding to subsets $Z_\a(\mu)$ for different $\mu.$
The determinant of this matrix is an entire function of $x,$
and asymptotically when $\<\a,x\> \to +\infty,\alpha\in R_+$ 
$$ \op{LHS } \sim \sum_{\nu\in Z_\a(\mu)}e^{n\<\a,x\>}
q^{-n\<\nu,\a\>}P_\nu(x)$$
Therefore, asymptotically the determinant of the matrix of this
system is equal to a Vandermonde-type
determinant, which is nonzero. It follows that the determinant
of the system of equations is nonzero for generic $x.$

Suppose $\mu$ is such that not all $P_\mu(x)$ is not identically zero
(i.e. we have a nontrivial solution of the system of equations).
Then for such $\mu$ we need to have
$\op{Card}\big( Z_\a(\mu) \big) > k_\a$ for all $\a.$
\enddemo

The rest is based on the following combinatorial Lemma.
\proclaim{Lemma 3.3}
If a subset $L\subset \Q_+$ is such that for any
$\a\in R,\mu\in L$ we have
$$\op{Card}\big( Z_\a(\mu) \cap L \big) > k_\a \tag 3-4$$
then any maximal weight $\mu_0\in L$ satisfies
$$\mu_0\ge\Theta = \sum_\a k_\a\cdot\a$$
\endproclaim

\demo{Proof of Lemma}
We start by proving that $\<\mu_0,\La_1\> \ge \<\Theta,\La_1\>.$
To any weight $\mu$ assign a string of numbers
$$\mu \mapsto \{\<\mu,\La_1\>,\dots,\<\mu,\La_r\> \}$$
where $\La_1,\La_2,\dots,\La_r$ are the fundamental weights.
Take the set of all roots $\beta\in R^+$ such that
$\<\beta,\La_1\> > 0$
and arrange them in the increasing lexicographical order:
$$\beta_1 \prec \beta_2 \prec \dots \prec \beta_N$$
Our goal is to construct a ``descending'' sequence of weights
$\mu_1,\dots,\mu_N \in L$ such that
$$\<\mu_{j-1},\La_1\> - \<\mu_j,\La_1\> \ge
k_{\beta_j}\<\beta_j,\La_1\> \tag 3-5$$
Then we conclude that
$$\<\mu_0,\La_1\> \ge \<\mu_N,\La_1\> +
\sum_{j=1}^N k_{\beta_j}\<\beta_j,\La_1\> \ge
0 + \sum_{\a\in R^+} k_\a\<\a,\La_1\> = \<\Theta,\La_1\> \tag 3-6$$

We now explain how to construct $\mu_{j+1}$ given $\mu_j.$

We will build an auxiliary sequence
$\mu_j^i,\quad i=0,1,2,\dots,k_{\beta_j}$
such that three conditions are satisfied:
$$\mu_j^i + n\beta_j + \sum_{l=2}^r n_l \a_l\in L, \  n,n_l\ge0,
\text{ implies all } n_l=0 \tag *$$
$$\mu_j^i + n \beta_j\in L
\text{ for at most } i \text{ different } n\in\N, \tag **$$
$$\<\mu_j^{i+1},\La_1\> \le \<\mu_j^i,\La_1\> - \<\beta_j,\La_1\> \tag ***$$

We take $\mu_j^0 = \mu_j.$
Suppose we have already constructed $\mu_j^i, i<k_{\beta_j},$
satisfying (*) and (**).
Conditions (3-4) and (**) guarantee that we can find in $L$ a weight
of the form $\mu_j^i - n\beta_j, n\in\N.$
Take the {\it minimal} positive integer $m$ such that the set
$$\{ \mu\in L | \ 
\<\mu,\La_1\> = \<\mu_j^i - m\beta_j,\La_1\>,\quad
 \mu \succeq \mu_j^i - m\beta_j \}$$
is nonempty, and let $\mu_j^{i+1}$ be its maximal element in the sense of
lexicographical order. We claim that $\mu_j^{i+1}$ satisfies
(*), (**) and (***). Indeed, (***) is clear from the construction.
To prove (*) we note that if we have
$$\mu = \mu_j^{i+1} + n \beta_j + \sum_{l=2}^r n_l \a_l\in L$$
for $n,n_l\ge 0$ such that not all of them are equal to zero,
then either $n\ge1$ and $m$ was not minimal, or $n=0$ and (*)
was not satisfied for $\mu_j^i$ - contradiction.
Finally, (**) is obvious when
$\mu_j^i - \mu_j^{i+1} = m \beta_j;$ if it is not the case, i.e.
$\mu_j^i - \mu_j^{i+1} = m \beta_j + \sum_l \a_l,$ then there is no
$n\in\N$ at all such that $\mu_j^{i+1} + n \beta_j\in L.$

Having thus constructed all
$\mu_j^i,\quad i=1,2,\dots,k_{\beta_j},$ we can put
$\mu_{j+1} = \mu_j^{k_{\beta_j}}.$
From construction we see that (3-5) is satisfied, and by (3-6)
we conclude that indeed $\<\mu_0,\La_1\> \ge \<\Theta,\La_1\>.$
We can prove similarly that
$$\<\mu_0,\La_i\> \ge \<\Theta,\La_i\>,\quad i=1,2,\dots,r$$
Therefore, $\mu_0 \ge \Theta.$ The Lemma is proven.\qed
\enddemo

\proclaim{Corollary 3.4}
The $\Psi$-function, satisfying (3-1), is unique up to a factor, depending on $x.$
\endproclaim
\demo{Proof}
If we have another function $\Psi'(\la,x)$ with highest coefficient
$P'_\Theta(x),$ satisfying (3-1), then the function 
$$\phi(\la,x) = \Psi'(\la,x) -
P'_\Theta(x) \(P_\Theta(x)\)^{-1} \Psi(\la,x)$$
will still satisfy (3-1), but its support will only contain
weights $\mu<\Theta.$ By Proposition 3.2, $\phi(\la,x)\equiv 0,$
and the statement follows. \qed
\enddemo

\proclaim{Proposition 3.5}
Any function $\phi(\la,x),$ satisfying (3-1), which has the form
$$\phi(\la,x) = e^{\<\la,x\>} Q(\la,x),$$
where $Q(\la,x)$ is an $\End$-valued $q$-polynomial with depending
on $x$ coefficients, can be represented as
$$\phi(\la,x) = c(\la,x)\Psi(\la,x)$$
for some $\End$-valued $q$-polynomial $c(\la,x)$
with depending on $x$ coefficients, satisfying
$$c\(\la-\frac{n\a}2,x\) = c\(\la+\frac{n\a}2,x\)$$
\endproclaim
The proof uses induction on support of $c(\la)$ and is very
similar to the one in [ES]. However, we will not need this
fact for our main theorem, and we leave the details to the reader.

\vskip .1in
\vskip .1in
\centerline{\bf 4. Existence of difference operators.}
\vskip .1in

Introduce a family of difference operators $T_\La,$ corresponding to
weights $\La\in \P,$ acting on $\h^*,$ by
$$T_\La (x) = x + \La \log q^2, \quad x\in\h^* $$
They naturally act on functions on $\h^*;$
for example, for $f(x) = e^{\<\a,x\>}$ we have
$$(T_\La f)(x) = e^{\<\a,x+\La \log q^2\>} = 
 q^{2\<\a,\La\>}e^{\<\a,x\>} = q^{2\<\a,\La\>}f(x)$$

\proclaim{Proposition 4.1}
Any function $\phi(\la,x),$ satisfying (3-1), which has the form
$$\phi = e^{\<\la,x\>} P(\la,x)$$
for some $q$-polynomial $P(\la,x),$ can be represented as
$$\phi(\la,x) = D \Psi(\la,x)$$
for a unique difference operator $D$ with depending on $x$ coefficients.
\endproclaim
\demo{Proof}
Uniqueness is obvious, since otherwise the $\Psi$-function would be
annihilated by a nontrivail difference operator for all $\la,$
which is impossible. (See, for example, [EK]).

To prove existence of $D$ we use induction on the support of $P(\la,x).$
Consider the family of all finite subsets $L \subset \Q$ such that
for $\mu,\nu\in \Q$
$$\mu\in L, \nu<\mu \Rightarrow \nu\in L$$
Such subsets are all finite unions of
cones of the form $\Cal C_\mu = \{ \nu\in\Q | \nu\le \mu \}.$

We prove that if our statement is true for all $q$-polynomials $P(\la,x)$
whose support is strictly contained in $L,$ then it is also true
for $q$-polynomials with support $L.$

If $L$ does not contain weights $\mu\ge\Theta$ then by Lemma
we have $P(\la,x) \equiv 0,$ and the operator $D\equiv 0$.

Suppose there is a $\nu\in L$ such that $\nu\ge\Theta.$
Consider the set of all such $\nu$'s and let $\mu$ be a maximal
element from this set.

Consider the function
$$\phi'(\la,x) = \phi(\la,x) - P_\mu(x) T_{\mu-\Theta} \Psi(\la,x)$$
It satisfies (3-1) and has support strictly contained in $L.$ By
induction hypothesis we can represent $\phi'$ as
$$\phi'(\la,x) = D' \Psi(\la,x)$$
The operator
$$D = P_\mu(x) T_{\mu-\Theta} + D'$$
satisfies the required properties. \qed
\enddemo

We now prove a simple technical
\proclaim{Lemma 4.2}
Let a $q$-polynomial $c(\la)$ be represented as
$$c(\la) = \sum_{\pi\in\P/\Q} c_\pi(\la)$$
where $c_\pi(\la)$ are $q$-polynomials with support
in the coset $\pi + \Q.$

Suppose for some $n,\a$ we have that
$$c\(\la+\frac{n\a}2\) = c\(\la-\frac{n\a}2\) \tag 4-1$$
whenever $\<\a,\la\>=0.$ Then (4-1) is also satisfied
for each $c_\pi(\la).$
\endproclaim
\demo{Proof}
Property (4-1) is equivalent to divisibility
by $q^{2\<\a,\la\>}-1$ of the $q$-polynomial
$$\tilde c(\la) = c\(\la+\frac{n\a}2\) - c\(\la-\frac{n\a}2\)$$
On the other hand, one can see that
$$\tilde c(\la) = \sum_{\pi\in\P/\Q} \tilde c_\pi(\la),$$
where
$$\tilde c_\pi(\la) = c_\pi\(\la+\frac{n\a}2\) -  c_\pi\(\la-\frac{n\a}2\)$$
Clearly, $\tilde c(\la)$ is divisible by $q^{2\<\a,\la\>}-1$
if and only if all $\tilde c_\pi(\la)$'s
are divisible by $q^{2\<\a,\la\>}-1,$ and the Lemma follows.
\qed\enddemo

\proclaim {Theorem 4.3}
For any $q$-polynomial $c(\la)$ such that for any $\a\in R^+$
$$c\(\la+\frac{n\a}2\) = c\(\la-\frac{n\a}2\),
   \quad n=1,2,\dots,k_\a \tag 4-2$$
whenever $\<\a,\la\>=0,$ there exists a difference operator $D_c$ with
coefficients in $\End$ such that
$$D_c \Psi(\la,x) = \Psi(\la,x) c(\la+\rho)$$
The correspondence $c(\la) \to D_c$ is a homomorphism of rings.
\endproclaim
{\bf Remark.} We put $c(\la+\rho)$ on the right since in that form
it admits generalization to the matrix case (see Theorem 4.4). Of
course, for scalar $q$-polynomial $c(\la)$ we could write it in a more
traditional form $D_c\Psi(\la,x) = c(\la+\rho)\Psi(\la,x).$
\demo {Proof}
By Lemma 4.2, it suffices to prove the theorem for $c_\pi(\la)$
of the form
$$c(\la) = q^{2\<\mu_0,\la\>} \sum_{\mu\in \Q_+}
c_\mu q^{-2\<\mu,\la\>}$$
for some $\mu_0\in\P.$ Consider the function
$$\phi(\la,x) = T_{-\mu_0}\Psi(\la,x) c(\la+\rho).$$
It satisfies (3-1), and it has the form
$$\phi(\la,x) = e^{\<\la,x\>} \sum_{\mu\in \Q_+}
q^{-\<\mu,\la\>} Q_\mu(x),$$
By Proposition 4.1 it can be represented as
$$\phi(\la,x) = D \Psi(\la,x)$$
for some difference operator $D.$ Put $D_c = T_{\mu_0} D.$
Then we have
$$\gather
D_c\Psi(\la,x) = T_{\mu_0} D \Psi(\la,x) = T_{\mu_0}\phi(\la,x) = \\
 T_{\mu_0} T_{-\mu_0}\Psi(\la,x) c(\la+\rho) = \Psi(\la,x) c(\la+\rho),
\endgather$$

We now prove the homomorphism property.
Suppose we have two polynomials $c(\la)$ and $c'(\la).$
It is easily checked that operator $D_{cc'} - D_c D_{c'}$
annihilates the $\Psi$-function for any $\la$; therefore it
has to be identically zero. Hence our correspondence is a
homomorphism of rings.
\qed\enddemo

We have big supply of (scalar) $q$-polynomials, satisfying (4-2),
arising from the algebra of Weyl group invariant $q$-polynomials,
which is freely generated by the Casimir elements $c_1,\dots,c_r.$
This give us $n$ algebraically independent difference
operators $D_1,\dots,D_r.$
However, there exist other $q$-polynomials, with this property.
For instance,
$$c_0(\la) = \prod_{\a\in R^+}\prod_{n=-k_\a}^{k_\a} 
\( q^{2\<\a,\la\> + n\<\a,\a\>} - 1 \)$$
also satisfies (4-2).
It gives rise to a difference operator, commuting with all those
generated by the Casimir elements, but not lying in the ring generated by them.
This procedure gives examples of what we called algebraically
integrable commutative rings of difference operators.

Theorem 4.3 can be slightly generalized to the matrix case.

\proclaim{Theorem 4.4}
\roster
For any $\End$-valued $q$-polynomial $C(\la)$ such that 
$$C\(\la+\frac{n\a}2\)B_{n\a}(\la) = B_{n\a}(\la)C\(\la-\frac{n\a}2\),
\qquad \a\in R^+, n=1,\dots,k_\a$$
whenever $\<\a,\la\>=0,$ there exists a unique difference operator $D_C$
with coefficients in $\End,$ such that
$$D_C \Phi(\la,x) = \Psi(\la,x) C(\la+\rho)$$
The correspondence $\xi: C(\la) \mapsto D_C$ is a homomorphism of rings.
\endroster
\endproclaim
\demo{Proof}
The argument used in proof of Theorem 4.3 in the obvious way extends
to the matrix case.
\qed\enddemo

{\bf Remark.} Operators $D_1,\dots,D_r$ act on $\Psi$-fuction as
scalars, and therefore commute with all operators $D_C$ constructed
as above. In fact, one can show that the centralizer of the subring
generated by operators $D_1,\dots,D_r$ in ${\frak D}_q^r(U[0])$
coincides with the image of $\xi.$
We do not include the proof of this statement here.

In the next section we explain how our
construction is related to Macdonald theory.

\vskip .1in
\vskip .1in

\centerline{\bf 5. Root system $A_n$ and Macdonald theory. }
\vskip .1in

Consider a special case of our construction for
$\g = \frak{sl}_N,\quad U=U_k=S^{(kN)}V$, where $V$ is the 
fundamental representation. 
It is well-known that the zero weight subspace $U[0]$ 
is one-dimensional, and we can regard $\Psi(\la,x)$
as a scalar-valued function.
Note also that in this case $k_\a = k$ for all $\a\in R^+,$
and $\Theta = k\sum_{\a\in R^+} \a = 2k\rho.$
We have $d_i=1,$ and therefore $q_i=q$ for all $i=1,\dots,r.$

We first prove an important property of partial traces
$$\Cal B_\mu(\la) = Tr \( \op{Proj}_{\ver[\la-\mu]} \circ \twin
 \circ \op{Proj}_{\ver[\la-\mu]} \),$$
introduced in Section 3. We will use the notation
$$[n] = \frac{q^n-q^{-n}}{q-q^{-1}}, \qquad
  [n]_+ = \frac{q^{2n}-1}{q^2-1}, \qquad
  [n]_- = \frac{1-q^{-2n}}{1-q^{-2}}$$

\proclaim{Proposition 5.1}
Given any $\a\in R^+, n=1,\dots,k,$ we have for all $\mu\in\Q$
$$\Cal B_{\mu}(\la) = q^{-\<n\a,\Theta\>} \cdot
 \Cal B_{\mu-n\a}(\la-n\a) \tag 5-1$$
whenever $\<\a,\la+\rho\>=n.$
\endproclaim
\proclaim{Corollary 5.2}
For $\a\in R^+, n=1,\dots,k,$ the function $\Psi(\la,x)$ satisfies
$$\Psi\(\la+\frac{n\a}2,x\) = q^{-n\<\a,\Theta\>} \Psi\(\la-\frac{n\a}2,x\)$$
whenever $\<\a,\la\> = n.$
\endproclaim

\demo{Proof of Proposition 5.1}
For given $n,\a$ it is sufficient to prove a special case
of (5-1), corresponding to $\mu = n\a:$
$$\Cal B_{n\a}(\la) = q^{-\<n\a,\Theta\>} \cdot \Cal B_0(\la-n\a)
 = q^{-\<n\a,\Theta\>} \chi(\la-n\a) \tag 5-2$$
Indeed, let $\la$ be such that $\<\a,\la+\rho\>=n.$ Then
the image of $\tilde\Phi^u_{\la}$ is contained in $M_{\la-n\a}\otimes U.$
From (5-2) we see that
$$\twin\bigr|_{M_{\la-n\a}} = q^{-\<n\a,\Theta\>} \cdot \tilde
\Phi_{\la-n\a}^u,$$
and the more general formula (5-1) follows.
We now use induction on the height of root $\a$ to prove
formula (5-2), and thus the Proposition 5.1.

{\bf Base of induction.}
Consider the case when $\a = \a_i$ is a simple root.
Let $\la$ be generic in the hyperplane $\<\a_i,\la+\rho\>=n.$
Then $\ver$ contains a unique nonzero proper submodule $\ver^1,$
generated by the singular vector $v_{\la-n\a_i} = F_i^n\hw.$
One can check that
$$\twin\hw =
 \frac {(q^{-1}-q)^n}{[n]_+!}
 \(\prod_{m=n+1}^k \chi_m^{\a_i}(\la)\)
 \(\prod_{\beta\ne\a_i} \prod_{m=1}^k  \chi_m^\beta(\la) \)
 v_{\la-n\a}\otimes E_i^n u + \text{t.o.w },$$
where $t.o.w.$ denotes ``terms of other weights'' (in the first component).

It follows that
$$\Cal B_{n\a_i}(\la) = \frac {(q^{-1}-q)^n}{[n]_+!}
 \(\prod_{m=n+1}^k (1-q^{-2(n-m)})\)
 \(\prod_{\beta\ne\a_i} \prod_{m=1}^k \chi_m^\beta(\la)\)
  F_i^n E_i^n$$

It is known that $F_i^n E_i^n$ acts as multiplication by
$\frac{[k+n]!}{[k-n]!}$ in $U[0].$ Since
$$\frac {(q^{-1}-q)^n}{[n]_+!} 
 \(\prod_{m=n+1}^k (1-q^{-2(n-m)})\)  \frac{[k+n]!}{[k-n]!} =
  q^{-2kn} \prod_{m=1}^k \(1 - q^{2(n+m)}\),$$
for $\la$ from the hyperplane $\<\a_i,\la+\rho\> = n$ we have
$$\gather
\Cal B_{n\a_i}(\la) = q^{-2kn} \prod_{m=1}^k \(1 - q^{2(n+k)}\)
 \(\prod_{\beta\ne\a_i} \prod_{m=1}^k \chi_m^\beta(\la)\) = \\
 q^{-2kn} \(\prod_{m=1}^k \chi_m^{\a_i}(\la-n\a_i)\)
 \(\prod_{\gamma\ne\a_i} \prod_{m=1}^k \chi_m^\gamma(\la-n\a_i)\) =
 q^{-n\<\a_i,\Theta\>} \chi(\la-n\a_i)
\endgather$$

{\bf Induction step.}
Suppose (5-1) is true for all roots $\beta$ such
that $\op{height}\beta < \op{height} \a.$ We are going to
prove that (5-2) is true also for $\a.$

We first show that
$B_{n\a}(\la)$ is divisible by factors $\chi_m^\beta(\la-n\a)$
for all $\beta\ne\a, m=1,\dots,k.$
It suffices to prove that $B_{n\a}(\la)$ vanishes whenever
$\<\beta,\la-n\a+\rho\>=m.$

Consider two cases:
\roster
\item
If $s_\a(\beta)\in R^+$, put $\gamma = s_\a(\beta).$ Then
$\<\gamma,\la+\rho\>=\<\beta,\la-n\a+\rho\>=m.$

If $\la$ is generic from hyperplane $\<\gamma,\la+\rho\>=m,$
then the image of $\tilde\Phi^u_\la$ is contained in $M_{\la-m\gamma}.$
Since $n\a-m\gamma\notin\Q_+$, there will be no terms contributing
to $\Cal B_{n\a}(\la) = Tr\bigr|_{\ver[\la-n\a]}\tilde\Phi^u_\la,$ and
$\Cal B_{n\a}(\la)=0$ generically (and, therefore, identically)
in the hyperplane $\<\gamma,\la+\rho\>=m.$

\item
If $s_\a(\beta)\notin R^+$, put
$\gamma = -s_\a(\beta)\in R^+.$ Then $\a=\beta+\gamma,$
and we can assume the induction hypothesis true for $\beta$ and
$\gamma.$

We have:
$$\<\gamma,\la+\rho\> = - \<\beta,\la-n\a+\rho\> = -m$$
$$\<\beta,\la+m\gamma+\rho\> = \<\beta,\la-n\a+\rho\> + n\<\beta,\a\>
 + m\<\gamma,\beta\> = m + n - m = n$$
Therefore,
$$\Cal B_{n\a}(\la) =
 q^{\<m\gamma,\Theta\>} \Cal B_{n\a+m\gamma}(\la+m\gamma) =
 q^{\<m\gamma-n\beta,\Theta\>}
 \Cal B_{(n+m)\gamma}(\la+m\gamma-n\beta)$$
Also,
$$\<\a,\la+m\gamma-n\beta+\rho\> = n + m - n = m$$
so the image of $\twin$ is contained in $M_{\la - (m+n)\beta},$ and
$$\Cal B_{(n+m)\gamma}(\la+m\gamma-n\beta) = 0$$
because $m\a - (m+n)\beta \notin \Q_+$ and there are no terms
contributing to $\Cal B_{(n+m)\gamma}(\la+m\gamma-n\beta).$
\endroster

We have proved that $\Cal B_{n\a}(\la)$ vanishes in the
required hyperplanes, and is therefore divisible
by all the required factors. Thus, in the hyperplane 
$\<\a,\la+\rho\> = n$ we get
$$\Cal B_{n\a}(\la) = C(\la)
\prod_{\beta\ne\a}\prod_{j=1}^k \chi_n^\beta(\la-n\a)$$
for some $q$-polynomial $C(\la).$ It is easy to see by
comparing highest terms that $C(\la)$ is constant in the hyperplane
$\<\a,\la+\rho\>=n.$ To compute this constant, take $\la$ generic
such that $\<\a,\la+\rho\>=\<\beta,\la+\rho\>=n.$
Then automatically $\<\gamma,\la+\rho-n\beta\> = n.$
We have:
$$\Psi(\la,x) = q^{-\<n\beta,\Theta\>}\Psi(\la-n\beta,x) =
q^{-\<n\beta,\Theta\>}q^{-\<n\gamma,\Theta\>}\Psi(\la-n\beta-n\gamma,x)$$
But $\la-n\beta-n\gamma = \la-n\a$ does not lie on any
Kac-Kazhdan hyperplanes, therefore $\Psi(\la-n\a,x)\ne0,$ and
$$C(\la) = q^{-\<n\beta,\Theta\>}q^{-\<n\gamma,\Theta\>} =
 q^{-n\<\a,\Theta\>}$$
Proposition 5.1 is now proved.
\qed\enddemo

Let $\la$ be a dominant integral weight, and
 $V_\la$ is the irreducible $\qgr$-module
with highest weight $\la$.

\proclaim{Proposition 5.3 (Generalized Weyl formula)}
The operator $\tilde\Phi^u_\la: M_\la\to M_\la\o U$ descends 
to a homomorphism $V_\la\to V_\la \o U$.
The function
$$\tilde p_\la(x) = Tr\bigr|_{V_\la}(\twin e^x) $$
 is expressed
in terms of the functions $\Psi(\la,x)$ by
$$q^{\<\Theta,\la\>}\tilde p_\la(x) = \sum_{w\in W} (-1)^{l(w)}
  q^{\<\Theta,w^\rho\la\>}  \Psi(w^\rho \la,x) \tag 5-3$$
\endproclaim
\demo{Proof}
The operator $\tilde\Phi^u_\la$ defines an operator
$M_\la\to V_\la\o U$. This operator has to factor through $V_\la$
because it lands in a finite dimensional representation. Thus, 
$\tilde\Phi^u_\la$ in fact defines an operator
$V_\la\to V_\la\o U$. 

Recall that for $\la\in\P_{++}$ we have a resolution
$$0 \gets V_\la\gets \ver^0 \gets \ver^1 \gets \ver^2 \gets \dots,$$
where
$$\ver^0 = \ver,\qquad \ver^i = \bigoplus_{l(w)=i} M_{w^\rho \la}$$
For matrix traces we have as for usual characters 
$$\tilde p_\la(x) = Tr\bigr|_{V_\la}(\twin e^x) =
\sum_i (-1)^i Tr\bigr|_{M_\la^i}(\twin e^x)$$

When $\la$ is generic from hyperplane
$\<\a_i,\la+\rho\> = \frac n2 \<\a_i,\a_i\>,\quad n>k,$ then
$\ver$ contains a submodule $\ver^1 \cong M_{\la-n\a_i},$
generated by a singular vector $v_{\la-n\a_i}=F_i^n\hw.$
Since $\ver[\la-n\a_i]$ is one-dimensional, we can write
$$\twin v_{\la-n\a_i} = v_{\la-n\a_i}\otimes u'' + \dots$$
Then we will have
$\twin\bigr|_{\ver^1} = \tilde\Phi^{u''}_{\la-n\a_i}.$
To compute $u''$ we use the formula
$$\twin\hw = \sum_{m=0}^k \frac{(q^{-1}-q)^m}{[m]_+!}
\( \prod_{l=m+1}^k\chi_l^{\a_i}(\la) \)
\( \prod_{\beta\ne\a_i}\prod_{l=1}^k \chi_l^\beta(\la) \)
F_i^m\hw \otimes E_i^m u + \dots,$$
where we retained only terms, which will
contribute to the expression for $u''.$
Applying $\Delta(F_i^n)$ to the RHS, and collecting terms
involving $v_{\la-n\a_i},$ we deduce that
$$\gather
u'' = \sum_{m=0}^k \qq{n \\ m} \frac{(q^{-1}-q)^m}{[m]_+!}
\( \prod_{l=m+1}^k \chi_l^{\a_i}(\la) \)
\( \prod_{\beta\ne\a_i}\prod_{l=1}^k \chi_l^\beta(\la) \)
  F_i^mK_i^{n-m}E_i^m u =\\
\sum_{m=0}^k (-1)^m q^{m-k}\qq{n \\ m} \frac{(q^{-1}-q)^k}{[m]_+!}
\frac{[n-m-1]_-!}{[n-k-1]_-!}
\( \prod_{\beta\ne\a_i}\prod_{l=1}^k \chi_l^\beta(\la) \)
 q^{2m(n-m)} F_i^mE_i^m u
\endgather$$
The operator $F_i^m E_i^m$ acts in $U[0]$ as multiplication by
$\frac{[k+m]!}{[k-m]!}.$ 

Also, we need to use the identity
$$ \sum_{m=0}^k (-1)^m \frac{[n]!}{[m]![n-m]!} 
 \frac{[n-m-1]!}{[m]![n-k-1]!}
\frac{[k+m]!}{[k-m]!} = (-1)^k \frac{[n+k]!}{[n]!}.$$
It can be interpreted as equality of two
polynomials in $z=q^{2n}$ of degree $k.$ To prove this identity it
suffices to check it for $n=0,\dots,k,$ when there is only one
nonzero term in the LHS of the equation.

After easy transformations, we conclude that
$$\gather
u'' = q^{-2kn} \( \prod_{m=1}^k (1-q^{-2(n+m)}) \)
\( \prod_{\beta\ne\a_i}\prod_{m=1}^k \chi_m^\beta(\la) \) u =\\
q^{-2kn} \( \prod_{m=1}^k \chi_m^{\a_i}(\la-n\a_i) \)
\( \prod_{\gamma\ne\a_i}\prod_{m=1}^k \chi_m^\gamma(\la-n\a_i) \) u =
q^{-n\<\a_i,\Theta\>} \chi(\la-n\a_i) u
\endgather$$
It follows that the restriction of $\twin$ to the submodule
$\ver^1 \cong M_{\la-n\a_i}$ coincides with 
$q^{-n\<\a_i,\Theta\>}\tilde\Phi^u_{\la-n\a_i}$
for $\la$ generic from hyperplane
 $\<\a_i,\la+\rho\> = \frac n2 \<\a_i,\a_i\>.$
Therefore, it is true for all $\la$ from that hyperplane.
Even more generally, for $\la\in\P_{++}$ we have 
$q^{\<\Theta,\la\>}\twin\bigr|_{M_{w^\rho\la}} = 
q^{\<\Theta,w^\rho\la\>}\tilde\Phi^u_{w^\rho\la}.$
Thus
$$q^{\<\Theta,\la\>}Tr\bigr|_{M^i}(\twin e^x) = \sum_{l(w)=i}
q^{\<\Theta,\la\>}Tr\bigr|_{M_{w^\rho \la}}(\twin e^x) =
\sum_{l(w)=i} q^{\<\Theta,w^\rho\la\>} \Psi\(w^\rho \la,x\)$$
Therefore
$$q^{\<\Theta,\la\>}\tilde p_\la(x) = \sum_{w\in W} (-1)^{l(w)}
q^{\<\Theta,w^\rho\la\>} \Psi\(w^\rho(\la),x\)$$
\qed
\enddemo

Introduce a normalized matrix trace by
$$\psi(\la,x) = \frac{q^{\<\Theta,\la-\rho\>}\Psi(\la-\rho,x)}
{\prod_{\a\in R^+}\prod_{i=1}^k
   (q^ie^{\<\a,x\>/2\>} - q^{-i}e^{-\<\a,x\>/2\>} )} $$
Condition (3-1) can be rewritten for the function $\psi(\la,x)$ as
$$\psi\(\la+\frac{n\a}2\) = \psi\(\la-\frac{n\a}2\), \tag 5-4$$
for all $\a\in R^+, \  n=1,\dots,k$ and for all
$\la$ such that $\<\a,\la\> = 0.$

\proclaim{Corollary (Macdonald polynomials)}
Combining (5-3) with results in [EK], we see that
(up to a factor) Macdonald polynomials are equal to
$$p_\la(x) = \sum_{w\in W} (-1)^{l(w)}
    \psi(w (\la + (k+1)\rho),x) \tag 5-5$$
\endproclaim

\proclaim{Proposition 5.4 (Macdonald operators)}
The function $\psi(\la,x)$ is the common eigenfunction for
Macdonald operators, corresponding to $t=q^{k+1}:$
$$\Cal M_i = \sum_{w\in W} w \( \prod_{\<\a,\La_i\>=1}
\frac {q^{k+1}e^{\<\a,x\>/2}-q^{-k-1}e^{-\<\a,x\>/2}}
{e^{\<\a,x\>/2}-e^{-\<\a,x\>/2}} T_{\La_i} \),\quad i=1,\dots,r $$
The corresponding eigenvalues are $W$-invariant $q$-polynomials
$c_i(\la) = \sum_{w\in W} q^{\<\la,w\La_i\>}:$
$$\Cal M_i \psi(\la,x) = c_i(\la) \psi(\la,x)$$
\endproclaim
\demo{Proof}
We already know that the $\Psi$-function is the common eigenfunction
for a family of commuting difference operators, corresponding to
$c(\la),$ satisfying (4-2). It is shown in [EK] by using central
elements in $\qgr,$ that for the $\Psi$-function normalized as above,
the operators, corresponding to elementary symmetric functions
$c_i(\la),$ are exactly Macdonald operators.
\qed\enddemo

Now we study relations between $\psi$-functions for
different values of $k.$ We use notation $\psi_k(\la,x)$
for the $\psi$-function constructed from representation
$U_k.$

\proclaim{Theorem 5.5 (Shift operators)} 
There exist difference operators $G_k$ such that
$$\psi_{k+1}(\la,x) = G_k \psi_k(\la,x)$$
These operators are $W$ invariant, and there action in the basis 
of Macdonald polynomials is given by the formula
$$
G_kp_{k,\l}(x)=p_{k+1,\l-\rho}(x), \l-\rho\in \bold P_+;\ G_kp_{k,\l}=0,
\l-\rho\notin P_+.
$$
\endproclaim
{\bf Remark.} Shift operators in the q-deformed case
were introduced by Cherednik,
\cite{Ch1}.

\demo{Proof}
The argument from the proof of Proposition 4.1 can be used to
prove that any function, satisfying (5-4), can be represented as
$$\phi(\la,x) = D \psi_k(\la,x)$$
for some difference operator $D.$ Applying this to the function
$\psi_{k+1}(\la,x),$ we prove the existence of an operator $G_k$
such that
$$\psi_{k+1}(\la,x) = G_k \psi_k(\la,x)$$

From the generalized Weyl formula  we get
$$\gather
G_k p_{k,\la} = \sum_{w\in W} (-1)^{l(w)}G_k\psi_{k}(w(\la+(k+1)\rho),x) = \\
\sum_{w\in W} (-1)^{l(w)} \psi_k(w(\la+(k+1)\rho),x) =
p_{k,\la-\rho}(x)
\endgather$$
If $\l-\rho$ is not dominant then the right hand side of the last formula is
zero by (5-4).

We see that $G_k$ maps Macdonald polynomials
to $W$-invariant functions. This implies that $G_k$ is
itself $W$-invariant.
\qed\enddemo

{\bf Remark.}
We saw that shift operators relate eigenfunctions for Macdonald
operators, corresponding to different (integral) values of
parameter $k.$ One can write it in the form
$$G_k \Cal M_i^{(k)} = \Cal M_i^{(k+1)} G_k \tag 5-6$$
One can check that $G_k$ analytically depends on $k,$ therefore
one can extend equality (5-6) to the case of arbitrary $k.$
This implies the existence of shift operators in the general case,
which is proven in \cite{Ch1} using representation theory of
double affine Hecke algebras.

Denote
$$\delta_k(x) = \prod_{\a\in R^+}\prod_{i=-k}^k
\(q^ie^{\<\a,x\>/2}-q^{-i}e^{-\<\a,x\>/2} \)$$

\proclaim{Theorem 5.6 (Duality)}
The function $\varphi_k(\la,x) = \delta_k(x) \psi_k(\la,x)$
is symmetric with respect to transformation
$q^{\<\a,\la\>} \leftrightarrow e^{\<\a,x\>}.$
\endproclaim
\demo{Proof}
The idea of the proof is the same as in \cite{VSC}. First
we prove that $\varphi_k$ as a function of $x$ satisfies condition
(4-2), and then duality will follow from the uniqueness property.

We already know that $\psi_k$ is the eigenfunction for the Macdonald
operators. Therefore, $\varphi_k(\la,x)$ is the eigenfunction for the
operators, obtained from Macdonald operators by conjugation by
$\delta_k(x).$ Such an operator, corresponding to a minuscule weight
$\La,$ is also $W$-invariant and has the form
$$\widetilde {\Cal M}_\La = \sum_{w\in W} f_\La(wx) T_{w(\La)},$$
where
$$f_\La(x) = \prod_{\<\beta,\La\>=1}
 \frac {q^{-k}e^{\<\beta,x\>/2}-q^k e^{-\<\beta,x\>/2}}
 {e^{\<\beta,x\>/2}-e^{-\<\beta,x\>/2}}$$

Fix a root $\a$ and denote $W_\a = \{w\in W | \<\a,w\La\> = 1\}.$ 
Let $x$ be such that $\<\a,x\> = 0.$
The $\phi$-function does not have singularities along that hyperplane.
Collecting all the singular terms in the equation
$$\widetilde {\Cal M}_\La \varphi_k(\la,x) =
 c_\La(\la) \varphi_k(\la,x),$$
which occur for $w\in W$ such that
$\<\a,w\La\> = \pm 1,$ we obtain
$$\sum_{w\in W_\a} f_\La(wx) T_{w(\La)}  \varphi_k(\la,x) + 
\sum_{w\in W_\a} f_\La(s_\a w x) T_{(s_\a w)(\La)} \varphi_k(\la,x) = 0
\tag 5-7$$
When $x$ belongs to hyperplane $\<\a,x\> = 0$ we have
$$w \(\prod_{\<\beta,\La\>=1} 
\frac {q^{-k-1}e^{\<\beta,x\>/2}-q^{k+1} e^{-\<\beta,x\>/2}}
 {e^{\<\beta,x\>/2}-e^{-\<\beta,x\>/2}} \) =
- (s_\a w) \(\prod_{\<\beta,\La\>=1} 
\frac {q^{-k}e^{\<\beta,x\>/2}-q^k e^{-\<\beta,x\>/2}}
 {e^{\<\beta,x\>/2}-e^{-\<\beta,x\>/2}} \),$$
therefore (5-7) simplifies to
$$\sum_{w\in W_\a} f_\La(wx)
 \(\varphi_k\(\la,x_w + \frac \a2\) - 
   \varphi_k\(\la,x_w - \frac \a2\) \) = 0$$
where
$$x_w = x + w(\La)^\perp = x + w(\La) - \frac \a2$$
also belongs to hyperplane $\<\a,x_w\> = 0.$
It follows that 
$$\varphi_k\(\la,x+\frac \a2\) - \varphi_k\(\la,x-\frac \a2\) = 0$$
identically when $\<\a,x\>=0,$ i.e. that condition (5-4)
is satisfied for $n=1.$ Taking $n$-th power of operators
$\widetilde {\Cal M}_\La,$ and repeating the same argument, we
can prove that it is satisfied for $n=1,\dots,k.$

From the obvious modification of the theorem about uniqueness of
the $\Psi$-function, we conclude that the function $\varphi_k(\la,x)$
transforms into itself when we interchange
$q^{\<\a,\la\>} \leftrightarrow e^{\<\a,x\>}.$
\qed\enddemo

{\bf Remark.} This duality result is closely related to 
the symmetry of the difference Fourier pairing defined recently
by Cherednik \cite{Ch2}.

\vskip .1in
\vskip .1in

\centerline{\bf Appendix. Example: $\g=\frak{sl}_2.$}
\vskip .1in

In this case we have only one root $\a,$ and we can make
identifications
$$e^{\<\la,x\>} \leftrightarrow e^{\la x}, \qquad
  e^{\<\a,x\>} \leftrightarrow e^{2x}, \qquad
  q^{\<\a,\la\>} \leftrightarrow q^{\la}$$
Operators $T,T^{-1}$ act on function $f(x)$ by
$$(Tf)(x) = f(x+1),\qquad (T^{-1}f)(x) = f(x-1)$$

\subhead Case $k=0$ \endsubhead
This is the simplest example, corresponding to the case 
$U$ - trivial representation. In that case
$$\twin\hw = \tw \hw = \hw\otimes u $$
$$\Psi(\la,x) = \frac{e^{\la x}}{1-e^{-2x}} =
    \frac{e^{(\la+1)x}}{e^x - e^{-x}}$$
$$\psi_0(\la,x) = \frac {e^{\la x}}{e^x - e^{-x}}$$

Formula (5-5) is the usual Weyl formula
$$\op{char} V_\la = \frac{ e^{(\la+1)x} - e^{-(\la+1)x}}
{e^x - e^{-x}}$$

The operator, corresponding to $W$-invariant polynomial
$c_1(\la) = q^{\la}+q^{-\la},$ is the
Macdonald operator
$$\Cal M_1 = \frac{qe^x - q^{-1}e^{-x}}{e^x - e^{-x}} T +
 \frac{qe^{-x} - q^{-1}e^x}{e^{-x} - e^x} T^{-1} $$
and we have
$$\Cal M_1\psi_1(\la,x) = c_1(\la)\psi_1(\la,x)$$

Condition (4-2) gives no restriction on $c(\la),$ and we can set
$c_0(\la) = q^{\la}.$ The corresponding operator $\Cal M_0,$
whose existence is predicted by Theorem 4.1, is equal to
$$\Cal M_0 = \frac{qe^x - q^{-1}e^{-x}}{e^x - e^{-x}} T,$$
and $\Cal M_1=\Cal M_0+\Cal M_0^{-1}$.

\subhead Case $k=1$ \endsubhead
Let $U$ now be the 3-dimensional representation. We have:
$$\tw \hw =
\hw\otimes u - \frac {q-q^{-1}}{1 - q^{-2\la}} F\hw\otimes Eu$$
$$\twin\hw = (1-q^{-2\la})\hw\otimes u - (q-q^{-1}) F\hw\otimes Eu$$
$$\gather
\Psi(\la,x) = \frac{e^{\la x}}{1-e^{-2x}} \( 1-q^{-2\la}
 - (q^2-q^{-2}) \frac {e^{-2x}}{1-q^{-2}e^{-2x}} \) =\\
 = \frac{e^{\la x}}{1-e^{-2x}}
 \( \frac {1-q^2e^{-2x}}{1-q^{-2}e^{-2x}} - q^{-2\la} \)
\endgather$$
$$\psi_1(\la,x) = \frac {e^{\la x}} {e^x - e^{-x}}
 \( \frac {q^{\la}}{qe^x - q^{-1}e^{-x}}
  - \frac {q^{-\la}} {q^{-1}e^x - qe^{-x}} \)$$
It is easy to check that indeed 
$$\Psi(0,x) = q^{-2} \Psi(-2,x),\qquad
\psi_1(1,x) =  \psi_1(-1,x)$$
Function $\psi_1(\la,x)$ is related to $\psi_0(\la,x)$ by
$$\psi_1(\la,x) = G_0 \psi_0(\la,x)$$
where the shift operator $G_0$ is equal to
$$G_0 = \frac 1{e^x - e^{-x}} \( T - T^{-1} \)$$

The Macdonald operator $\Cal M_1$ is equal to
$$\Cal M_1 = \frac{q^2e^x - q^{-2}e^{-x}}{e^x - e^{-x}} T +
 \frac{q^2e^{-x} - q^{-2}e^x}{e^{-x} - e^x} T^{-1} $$

$$\Cal M_1\psi_1(\la,x) = c_1(\la)\psi_1(\la,x)$$
The operator $\Cal M_0,$ corresponding to the eigenvalue 
$c_0(\la)=q^{3\la}-[3]_qq^\la,$ is equal to
$$\Cal M_0 = \frac{q^4e^x - q^{-4}e^{-x}} {qe^x - q^{-1}e^{-x}}
  \frac {q^3e^x - q^{-3}e^{-x}} {e^x - e^{-x}} T^3 - 
[3]_q \frac {qe^x - q^{-1}e^{-x}} {q^{-1}e^x - qe^{-x}} T$$

\subhead Classical limit $q\to 1$ \endsubhead
Set $\epsilon = \log(q).$ We have the expansion
$$\psi_1 (\la,x) = 2\epsilon \cdot\psi_1^{(0)}(\la,x) +
 O(\epsilon^2) $$
where
$$\psi_1^{(0)}(\la,x) = \frac{e^{\la x}}{\(e^x - e^{-x}\)^2}
 \(\la - \frac{e^x + e^{-x}}{e^x - e^{-x}}\)$$
is the $\psi$-function for the classical case (cf. [ES]).
The difference operators become
$$\Cal M_1 = 2 + \epsilon^2 (\Cal D_2 + 4) + O(\epsilon^3)$$
$$\Cal M_0 = -2 + \epsilon^2 (3\Cal D_2 + 8) +
4\epsilon^3 (\Cal D_3 + 1) + O(\epsilon^4),$$
where commuting differential operators $\Cal D_2, \Cal D_3$ are equal
$$\Cal D_2 = \frac {\partial^2}{\partial x^2} + 
4 \frac {e^x + e^{-x}}{e^x - e^{-x}} \frac \partial {\partial x}$$
$$\Cal D_3 =  \frac {\partial^3}{\partial x^3} +
6 \(\frac {e^x + e^{-x}}{e^x - e^{-x}}\) \frac{\partial^2}{\partial x^2} +
 \(11 + \frac {12}{(e^x - e^{-x})^2}\) \frac\partial{\partial x} +
6 \( \frac{e^{3x}-3e^x-3e^{-x}+e^{-3x}} {(e^x-e^{-x})^3} \)$$

\vskip .1in

\vskip .1in

\end
\vskip .1in

\Refs

\ref \by [Ch1] I.Cherednik
\paper Double affine Hecke algebras and Macdonald's conjectures
\jour Annals of Math\vol 141\pages 191-216\yr 1995\endref

\ref\by [Ch2] I.Cherednik
\paper Macdonald's evaluation conjectures and difference Fourier transform
\jour Inv. Math.\vol 122\pages 119-145\yr 1995\endref

\ref\by [CV1] O.Chalykh, A.Veselov
\paper Commutative rings of partial differential operators and Lie algebras
\jour Comm.Math.Phys  \vol 126 \yr 1990
\pages 597-611 \endref

\ref \by [CV2] O.Chalykh, A.Veselov
\paper Integrability in the theory of Schrodinger operator and 
harmonic analysis
\jour Comm.Math.Phys  \vol 152 \yr 1993
\pages 29-40 \endref

\ref \by [CK] C.deConcini, V.Kac \paper Representations of quantum groups 
at roots of 1, 
\jour Operator algebras, unitary representations, enveloping algebras, 
and invariant theory\publ Birkhauser\publaddr Boston \pages 471-506
\yr 1990
\endref

\ref \by [EK] P.Etingof, A.Kirillov, Jr.
\paper Macdonald polynomials and representations of quantum groups
\jour Nucl Phys B \vol 317 \yr 1994
\pages 215-236 \endref

\ref \by [ES] P.Etingof, K.Styrkas
\paper Algebraic integrability of Schr\"odinger operators and
representations of Lie algebras
\jour Compositio Mathematica \vol 98\pages 91-112\yr 1995
\endref

\ref \by {FV} G.Felder and A.Varchenko \paper
Three formulas for eigenfunctions of integrable Schr\''odinger operators
\jour Preprint\yr 1995\endref

\ref \by [VSC] A.Veselov, K.Styrkas, O.Chalykh
\paper Algebraic integrability for the Schr\"odinger equation
and finite reflection groups
\jour Theor. and Math. Physics \vol 94 \issue 2 \yr 1993
\endref
